\documentclass[journal,transmag]{IEEEtran}
\ifCLASSINFOpdf
\else
\fi
\usepackage{hyperref,color,multirow}
\usepackage[top=0.7in, left=0.65in]{geometry}
\usepackage{graphicx}

\usepackage{amsmath,amsfonts,amssymb,amscd,bm,bbm}

\usepackage{multicol}
\usepackage{blindtext}

\newcommand{\opgrad}{\operatorname{grad}}
\newcommand{\opcurl}{\operatorname{curl}}
\newcommand{\opdiv}{\operatorname{div}}
\newcommand{\MR}{\mathbbm R}  

\begin{document}
%
\title{Multiscale Finite Element Formulations \\ for 2D/1D Problems}



\author{\IEEEauthorblockN{Karl Hollaus\IEEEauthorrefmark{1} and Markus Sch\"{o}binger\IEEEauthorrefmark{1}}
\IEEEauthorblockA{\IEEEauthorrefmark{1}Technische Universit\"{a}t Wien, Institute for Analysis and Scientific Computing, A-1040 Vienna, Austria} 
\thanks{Manuscript received Oct 26, 2022; revised xxx, xx. 
Corresponding author: K. Hollaus (email: karl.hollaus@tuwien.ac.at)}}

\markboth{}%
{Shell \MakeLowercase{\textit{et al.}}: Bare Demo of IEEEtran.cls for IEEE Transactions on Magnetics Journals}
%



\IEEEtitleabstractindextext{%
\begin{abstract}
Multiscale finite element methods for 2D/1D problems have been studied in this work to demonstrate their excellent ability to solve real-world problems. These methods are much more efficient than conventional 3D finite element methods and just as accurate. The 2D/1D multiscale finite element methods are based on a magnetic vector potential or a current vector potential. Known currents for excitation can be replaced by the Biot-Savart-field. Boundary conditions allow to integrate planes of symmetry. All presented approaches consider eddy currents, an insulation layer and preserve the edge effect. A segment of a fictitious electrical machine has been studied to demonstrate all above options, the accuracy and the low computational costs of the 2D/1D multiscale finite element methods.
\end{abstract}

\begin{IEEEkeywords}
Biot-Savart-field, eddy currents, edge effect, thin iron sheets, 2D/1D multiscale finite element method MSFEM
\end{IEEEkeywords}
}

\maketitle

\IEEEdisplaynontitleabstractindextext

%
\IEEEpeerreviewmaketitle

\section{Introduction}
\IEEEPARstart{T}{he} overall dimensions of a laminated core of electrical machines are essentially larger than the thickness of a single iron sheet and thus such machines represent a multiscale problem. Neglecting the magnetic stray fields at the end region, all iron sheets are exposed to roughly the same field distribution. Thus, a simulation of a single sheet instead of the whole core suffices. \\
The brute force way is to exploit three-dimensional (3D) finite element methods (FEMs) for the single sheet. However, simulations with 3D FEMs may still become too expensive for routine tasks \cite{HandBiroBelaDlal:13}. \\
Attractive alternative options to 3D FEMs are space splitting two-dimensional/one-dimensional (2D/1D) methods, see for instance \cite{BottChiamp:02}, \cite{PippArkk:09} and \cite{HenSteHamGeu:15}. Some of them only solve a static magnetic field problem or suffer from a high number of subdivisions along the thickness of the sheet to properly represent the penetration depth. The edge effect (EE) is ignored. The EE is particularly important, for instance, in the tooth tips which are exposed to high flux variations (\cite{HandBiroBelaDlal:13}) and because of the degrading effect due to the cutting process of iron sheets, see for example \cite{VandenJacobsHenrotteHameyer:2010}, \cite{BaliMuetze:17} and \cite{SundNairLehikArkkBelah:20}. \\ 
\begin{figure}
		\begin{minipage}{0.99\columnwidth}
			\begin{minipage}{0.65\columnwidth}
					\includegraphics[scale=0.32]{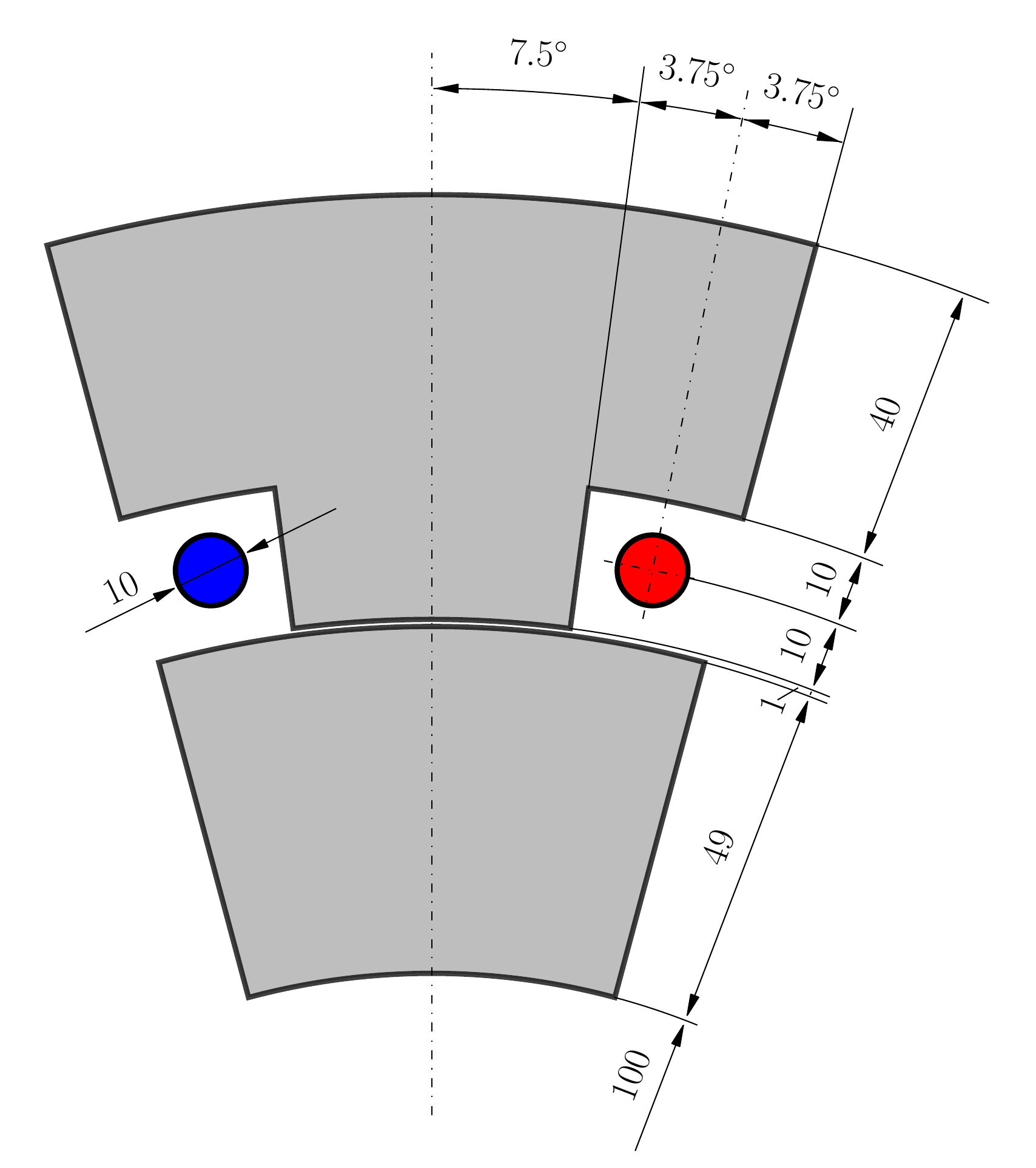} 
					\includegraphics[scale=0.15]{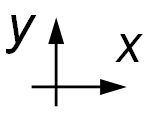} 
			\end{minipage}	
			\begin{minipage}{0.33\columnwidth}
					\includegraphics[scale=0.21]{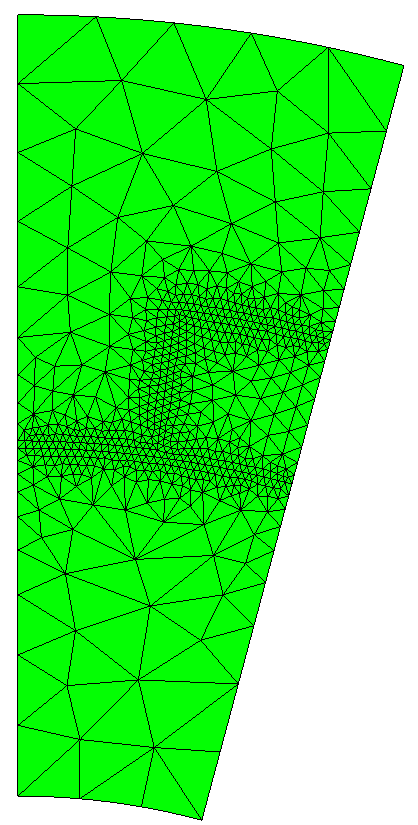} 
			\end{minipage}
		\end{minipage}
\caption{The segment is a twelfth of a fictitious electrical machine (left) with a rotor and a stator separated by an air gap and assumed to be in the xy-plane. Dimensions are in mm. Thickness of the iron sheets is $d$=0.5mm. Known opposed currents in conductors with circular cross-section represent the excitation. Finite element mesh of one half of the segment (right).}
\label{MotorSegmentGeo}
\end{figure}
A very efficient approach is the use of an effective material with a complex-valued magnetization curve \cite{SchoebTsukHoll:21}. However, methods with an effective material are restricted to problems in the steady state. \\
Therefore, the idea is to replace the current 3D FEMs or the 2D/1D methods by 2D/1D MSFEMs based on a magnetic vector potential (MVP) $\boldsymbol{A}$ or a current vector potential (CVP) $\boldsymbol{T}$ formulation. 
The 2D/1D MSFEM approaches require much fewer unknowns and the system matrix is much sparser than for 3D FEMs, which is very important for drastically reducing the computational cost. They even reduce the effort of the former 2D/1D methods significantly, while being as versatile as possible. \\ 
A 2D/1D MSFEM using trigonometric functions across the thickness of the sheet can be found in \cite{RasiloDlaPippBelArkk:11}. 
Our methods also have to consider eddy currents including the EE, see \cite{SchoeSchoeHoll_2D1D:19} and \cite{HollSchoebi2D1DEE:20}, account for an insulation layer in between the iron sheets, facilitate boundary conditions (BCs) to exploit planes of symmetry, and use Biot-Savart-fields (BSF) to avoid modeling of conductors carrying known currents. So far, an excitation has been introduced only by proper BCs in the tiny problem in \cite{SchoeSchoeHoll_2D1D:19} and \cite{HollSchoebi2D1DEE:20}. \\
First, the basic eddy current problem (ECP) with BCs is presented in Sec.~\ref{EddyCurrentProblem}. The segment of a fictitious electric machine in Fig.~\ref{MotorSegmentGeo} serves as model problem. For the sake of simplicity linear material relations and steady state are assumed, thus the work is carried out in the frequency domain. Nevertheless all advantages of 2D/1D MSFEMs over 3D FEM can be shown.
Then, four 2D/1D MSFEM approaches are presented in Sec.~\ref{2D1D_MSFEM_Forms}. In contrast to \cite{SchoeSchoeHoll_2D1D:19} and \cite{HollSchoebi2D1DEE:20}, the new approaches (\ref{eqT1}) and (\ref{eqA2}) use either $H(\opcurl)$ or $H^1$ finite element spaces. To evaluate the accuracy and efficiency of the  2D/1D MSFEMs, mixed FEMs have been used which are briefly discribed in Sec.~\ref{RefSols}. Simulation results obtained by the 2D/1D MSFEMs by means of the numerical example in Fig.~\ref{MotorSegmentGeo} are presented in Sec.~\ref{NumSimulations}. The accuracy of the 2D/1D MSFEMs by means of eddy current losses, the required unknowns, the non-zero entries in the system matrix and the computation times are shown. The accuracy of modeling the EE by $\boldsymbol{A}$ and $\boldsymbol{T}$ formulations has also been studied. In summary, the 2D/1D MSFEMs show high accuracy, comparable to the expensive 3D FEMs, but require very low computational cost.
\section{Eddy Current Problem} \label{EddyCurrentProblem}
\begin{figure}[ht]
	\begin{minipage}{0.99\columnwidth}
		\begin{center}
			\begin{minipage}{0.49\columnwidth}
				\includegraphics[scale=0.305]{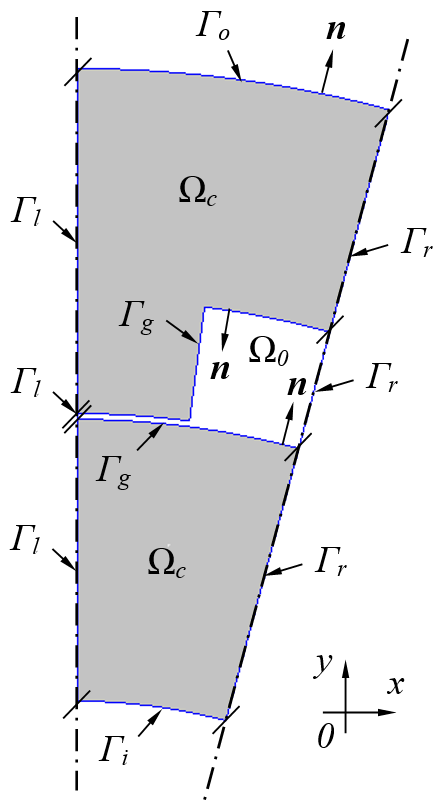} 
			\end{minipage}	
			\begin{minipage}{0.40\columnwidth}
				\begin{center}
					\includegraphics[scale=0.305]{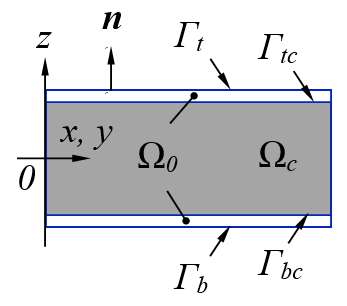} 
				\end{center} 
			\end{minipage}	
		\end{center}
	\end{minipage}
	\caption{One half of the segment in the $xy$-plane (left), iron sheet in grey. Detail of the cross-section with iron sheet and with half a layer of insulation on top and half on the bottom (right), not to scale.}
	\label{MotorSegmentBVP}
\end{figure}

An ECP has to be solved, see Fig.~\ref{MotorSegmentGeo}. The entire domain $\Omega=\Omega_{c} \cup \Omega_{0} \in \MR^3$ consists of the conducting domain (iron sheets) $\Omega_{c}$ and the nonconducting domain (air or insulation layer) $\Omega_{0}$, compare with Fig.~\ref{MotorSegmentBVP}. The normal vector $\boldsymbol{n}$ points out of $\Omega$ and $\Omega_c$, respectively. To facilitate the representation of the ECP the following definitions are introduced
\begin{eqnarray}
	\Gamma_H &=& \Gamma_r \label{BC_H} \\
	\Gamma_{J} &=& \Gamma_{ic} \cup \Gamma_{rc} \cup \Gamma_{gc} \cup \Gamma_{oc} \cup \Gamma_{bc} \cup \Gamma_{tc}, \label{BC_J} \\
	\Gamma_{B} &=& \Gamma_{i} \cup \Gamma_{o} \cup \Gamma_{l} \cup \Gamma_{b} \cup \Gamma_{t} \text{ and} \label{BC_B} \\
	\Gamma_E &=& \Gamma_{lc}. \label{BC_E}
\end{eqnarray}
The additional index $c$ means the part of the respective boundary connected with a conductive domain. 
The ECP consists of a quasi-static magnetic field
\begin{align}
	\opcurl \boldsymbol{H} = & \ \boldsymbol{J}, \label{EQ_1} \\
	\opcurl \boldsymbol{E} = & \ - j \omega \boldsymbol{B} \text{ and} \label{EQ_2} \\
	\opdiv \boldsymbol{B} = & \ 0 \label{EQ_3}
\end{align}
with
\begin{align}
	\boldsymbol{J} = \sigma \boldsymbol{E} \text{ or } \boldsymbol{E} = \rho \boldsymbol{J} \text{ and } \boldsymbol{B} = \mu \boldsymbol{H}  \text{ or } \boldsymbol{H} = \nu \boldsymbol{B}  \label{EQ_4}
\end{align}
in $\Omega_c$ and a static magnetic field
\begin{align}
	\opcurl \boldsymbol{H} = & \  \boldsymbol{J}_0 \text{ and} \label{EQ_5} \\
	\opdiv \boldsymbol{B} = & \ 0 \label{EQ_6}
\end{align}
with
\begin{align}
	\boldsymbol{B} = \mu \boldsymbol{H}  \text{ or } \boldsymbol{H} = \nu \boldsymbol{B} \label{EQ_7}
\end{align}
in $\Omega_0$ and the BCs
\begin{eqnarray}
	\boldsymbol{H} \times \boldsymbol{n} = \boldsymbol{0}& \ &\text{on } \ \Gamma_{H}, \label{EQ_8} \\
	\boldsymbol{J} \cdot \boldsymbol{n} = 0& \ &\text{on } \ \Gamma_{J},  \label{EQ_9} \\
	\boldsymbol{B} \cdot \boldsymbol{n} = 0& \ &\text{on } \ \Gamma_{B} \text{ and } \label{EQ_10} \\
	\boldsymbol{E} \times \boldsymbol{n} = \boldsymbol{0}& \ &\text{on } \ \Gamma_{E}, \label{EQ_11} \\
	\end{eqnarray}
where $\boldsymbol{H}$ is the magnetic field strength, $\boldsymbol{J}$ the electric current density, $\boldsymbol{E}$ the electric field strength, $\boldsymbol{B}$ the magnetic flux density, $\sigma$ the electric conductivity, $\rho$ the electric resistivity, $\mu$ the magnetic permeability, $\nu$ the magnetic reluctivity and $\boldsymbol{J}_0$ the prescribed electric current density, respectively, $j$ denotes the imaginary unit and $\omega$ the angular frequency. Since there is only the magnetic field in the entire domain $\Omega$, the continuity conditions on the interface $\Gamma_{c0}$ between $\Omega_{c}$ and $\Omega_{0}$ are
\begin{equation}
	\boldsymbol{H} \times \boldsymbol{n} \text{ and } \boldsymbol{B} \cdot \boldsymbol{n} \label{EQ_12}
\end{equation}
which are continuous. \\
If the symmetry in the segment is not used, $\Gamma_{l}$ does not exist. In the considered example, the BC $\boldsymbol{B} \cdot \boldsymbol{n} = 0$ on $\Gamma_{i} \cup \Gamma_{o}$ represents a simplification and is in reality only approximately fulfilled and on $\Gamma_{b} \cup \Gamma_{t}$ an assumption that no magnetic stray field exists.
\section{2D/1D MSFEM Formulations} \label{2D1D_MSFEM_Forms}
The entire domain $\Omega=\Omega_m \cup \Omega_0 \in \MR^2$ is composed of a laminated domain $\Omega_m$ representing an iron sheet and half an insulation layer on each side of the sheet, compare with the detail in Fig.~\ref{MotorSegmentBVP}, and the non-conducting domain $\Omega_0$ representing the air gap between the rotor and the stator and the space for conductors with known currents, see also Fig.~\ref{MotorSegmentGeo}. Note that the meaning of $\Omega_m$ and $\Omega_0$ depends on the context, either that of a 2D/1D MSFEM or that of the reference solution. 
\subsection{2D/1D MSFEM Approaches} \label{2D1D_MSFEM_Approaches}
\begin{figure}[ht]
\begin{center} \includegraphics[scale=0.20]{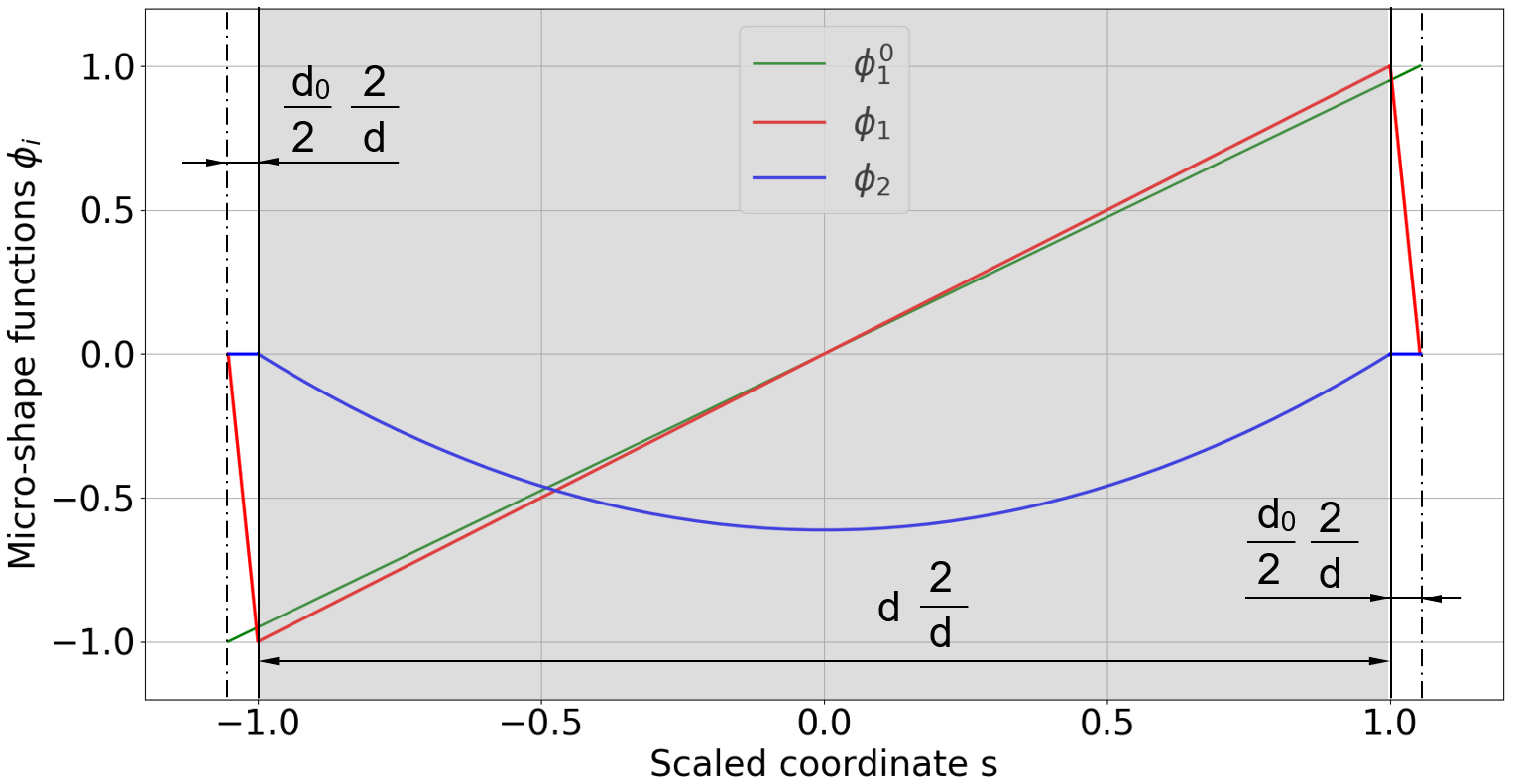} \end{center} 
\caption{Micro-shape functions $\phi_i$: The gray interval [-1, 1] represents the iron sheet, and
beyond that up to the dashed–dotted line, there is the insulation layer.}
\label{ShapeFunctions_phi10_ph1_ph2}
\end{figure}

The considered MSFEMs approaches 
%
\begin{align}
		\widetilde{\boldsymbol{T}}_{1} &= \boldsymbol{T}_0 + \phi_2 \boldsymbol{T}_2 + \boldsymbol{H}_{BS} \label{eqT1} \\
		\widetilde{\boldsymbol{T}}_{2} &= \opgrad \Phi_0 + \phi_2 \boldsymbol{T}_2 + \boldsymbol{H}_{BS} \label{eqT2} \\
		\widetilde{\boldsymbol{A}}_{1} &= \phi_1^0 \opgrad u_{10} + \phi_1 \boldsymbol{A}_1 + \opgrad (\phi_1 w_1) \label{eqA1} \\ 
		\widetilde{\boldsymbol{A}}_{2} &= \phi_1^0 \opgrad u_{10} + \phi_1 \opgrad u_{1} + \phi_{1,z} (0,0,w_1)^T  \label{eqA2}
\end{align}
are denoted by a tilde, where $\boldsymbol{H}_{BS}$ stands for the Biot-Savart field. The micro-shape functions (MSFs) $\phi_1$, $\phi_2$ and $\phi_1^0$, where $\phi_{1,z}$ is the derivative of $\phi_{1}$ with respect to $z$ are shown in Fig.~\ref{ShapeFunctions_phi10_ph1_ph2}.
 \\
Figure~\ref{ShapeFunctions_phi10_ph1_ph2} shows how the MSFs fit into the periodic structure with $p=d+d_0$, where $d$ is the thickness of the iron sheet and $d_0$ that of the insulation layer. The polynomials
\begin{align}
		\phi_1(s) &= s, \ \ \phi_2(s) = \frac{1}{2}\sqrt{\frac{3}{2}}(s^2-1) \label{eq:shapefunctions} 
\end{align}
are used as MSFs with the mapping $s=2z/d$, where $s \in [-1,1]$ and $z \in [-d/2,d/2]$. The MSFs $\phi_1^0$ and $\phi_1$ are extended linearly and become $\pm 1$ and $0$ in $\{-(d+d_0)/2,(d+d_0)/2\}$, respectively, and $\phi_2$ is extended by zero in $[-(d+d_0)/2,-d/2)$ and $(d/2,(d+d_0)/2]$ which includes the insulation layer. These polynomials facilitate the required tangential continuity of the unknowns in the multiscale approaches and $\phi_1^0$ allows to prescribe essential BCs. The required symmetry of the solution with respect to $z=0$, is ensured by selecting either even or odd MSFs in the 2D/1D MSFEM approaches explicitly. \\
\subsection{Properties of the MSFEM Formulations} \label{MSFEMProperties}
In order to facilitate the discussion of the 2D/1D MSFEM approaches and to find the true BCs more easily some intermediate results are presented:
\begin{align}
		\opcurl \boldsymbol{T}_{0} &= \left(\begin{array}{c} 0 \\ 0 \\ T_{0y,x}-T_{0x,y} \end{array}\right) \label{eq2D1Daux_1} \\
		\opcurl  (\phi_2\boldsymbol{T}_{2}) &= \left(\begin{array}{c} -\phi_{2,z}T_{2y} \\ \phi_{2,z}T_{2x} \\ \phi_{2}(T_{2y,x}-T_{2x,y}) \end{array}\right) \label{eq2D1Daux_3} \\
		\opcurl \boldsymbol{H}_{BS} &= \boldsymbol{J}_{0} \label{eq2D1Daux_4} \\
		\opcurl  (\phi_1^0 \opgrad u_{10}) &= \left(\begin{array}{c} -\phi_{1,z}^0 u_{10,y} \\ \phi_{1,z}^0 u_{10,x} \\ 0 \end{array}\right) \label{eq2D1Daux_5} \\
		\opcurl  (\phi_1\boldsymbol{A}_{1}) &= \left(\begin{array}{c} -\phi_{1,z}A_{1y} \\ \phi_{1,z}A_{1x} \\ \phi_{1}(A_{1y,x}-A_{1x,y}) \end{array}\right) \label{eq2D1Daux_6} \\
		\opcurl  (\phi_1 \opgrad u_{1}) &= \left(\begin{array}{c} -\phi_{1,z} u_{1,y} \\ \phi_{1,z} u_{1,x} \\ 0 \end{array}\right) \label{eq2D1Daux_7} \\
		\opcurl  (0,0,\phi_{1,z} w_1) &= \left(\begin{array}{c} \phi_{1,z} w_{1,y} \\ -\phi_{1,z} w_{1,x} \\ 0 \end{array}\right) \label{eq2D1Daux_8}
\end{align}
Magnetic fields that occur in air are represented by $\boldsymbol{T}_{0}$ in (\ref{eqT1}). As can be seen in (\ref{eq2D1Daux_3}) the term $\phi_2\boldsymbol{T}_{2}$ is essential for the $\boldsymbol{T}$ formulations (\ref{eqT1}) and (\ref{eqT2}) to get proper eddy current density distributions. The trace of $\boldsymbol{T}_{2}$ is zero and $J_z=\phi_{2}(T_{2y,x}-T_{2x,y})$ for the EE. This is a big advantage over the $\boldsymbol{A}$ formulations (\ref{eqA1}) and (\ref{eqA2}), which additionally require a third term. To preserve the EE, $\phi_2 \boldsymbol{T}_2$ can not be replaced by a term like $\phi_2 \opgrad (u_2)$, because $\opcurl(\phi_2 \opgrad (u_2))$ does not have a $z$-component. The EE is represented by $\phi_1 \opgrad u_{1}$ and by $\phi_{1,z} (0,0,w_1)^T$  in (\ref{eqA1}) and (\ref{eqA2}), respectively. Laminar currents due to $\phi_1^0 \opgrad u_{10}$ generate a total magnetic flux (\ref{eq2D1Daux_5}), which is perturbed either by (\ref{eq2D1Daux_6}), see \cite{HollSchoebi2D1DEE:20}, or by (\ref{eq2D1Daux_7}). The $z$-component of $\opcurl  (\phi_1\boldsymbol{A}_{1})$ provides a smoothing of the magnetic field at the transition from the iron sheet to the air. While the $z$-component of the current density due to $\phi_{1,z} w_1$ is accompanied by a magnetic field (\ref{eq2D1Daux_8}), $\opgrad (\phi_1 w_1)$ does not yield a magnetic field. 
\subsection{Boundary Conditions of the 2D/1D MSFEMs} \label{BCs_2D1D_MSFEMs}
\subsubsection{TMS1} \label{TMS1_BVPs}
\begin{eqnarray}
	\boldsymbol{T}_0 \times \boldsymbol{n} = \boldsymbol{0}&&\text{on } \Gamma_H \label{TMS1_BC_1} \\
	\boldsymbol{T}_2 \times \boldsymbol{n} = \boldsymbol{0}&&\text{on } \Gamma_{J_{2D}}, \label{TMS1_BC_2}  \\
	-j\omega\mu \widetilde{\boldsymbol{T}}_{1} \cdot \boldsymbol{n} = 0&&\text{on } \Gamma_{B_{2D}}, \label{TMS1_BC_3}  \\
	\opcurl  (\boldsymbol{T}_{0}+\Phi_2\boldsymbol{T}_{2}) \times \boldsymbol{n}= \boldsymbol{0}&&\text{on } \Gamma_E \label{TMS1_BC_4} 
\end{eqnarray}
\subsubsection{TMS2} \label{TMS2_BVPs}
\begin{eqnarray}
	\Phi_0 = 0&&\text{on } \Gamma_H \label{TMS2_BC_1} \\
	\boldsymbol{T}_2 \times \boldsymbol{n} = \boldsymbol{0}&&\text{on } \Gamma_{J_{2D}} \label{TMS2_BC_2}  \\
	-j\omega\mu \widetilde{\boldsymbol{T}}_{2} \cdot \boldsymbol{n} = 0&&\text{on } \Gamma_{B_{2D}} \label{TMS2_BC_3}  \\
	\opcurl  (\Phi_2\boldsymbol{T}_{2}) \times \boldsymbol{n}= \boldsymbol{0}&&\text{on } \Gamma_E \label{TMS2_BC_4} 
\end{eqnarray}
\subsubsection{AMS1} \label{AMS1_BVPs}
\begin{eqnarray}
	u_{10} = 0,&& \nonumber \\
	\boldsymbol{A}_1 \times \boldsymbol{n} = \boldsymbol{0}&&\text{on } \ \Gamma_{H} \label{AMS1_BC_3} \\
	-j\omega\sigma \widetilde{\boldsymbol{A}}_{1} \cdot \boldsymbol{n} = 0&&\text{on } \Gamma_{J_{2D}} \label{AMS1_BC_1}  \\
	\boldsymbol{A}_1 \times \boldsymbol{n} = \boldsymbol{0}&&\text{on } \ \Gamma_{B} \label{AMS1_BC_2} \\
\end{eqnarray}
\subsubsection{AMS2} \label{AMS2_BVPs}
\begin{eqnarray}
	-j\omega\sigma \widetilde{\boldsymbol{A}}_{2} \cdot \boldsymbol{n} = 0&&\text{on } \Gamma_{J_{2D}} \label{AMS2_BC_1}  \\
	u_{10} = 0&& \nonumber \\
	u_{1} = 0&& \nonumber \\
	w_{1} = 0&&\text{on } \ \Gamma_{B} \label{AMS2_BC_2} \\
	u_{10} = 0&& \nonumber \\
	u_{1} = 0&& \nonumber \\
	w_{1} = 0&&\text{on } \ \Gamma_{H} \label{AMS2_BC_3}
\end{eqnarray}
To obtain the respective weak form of the 2D/1D MSFEMs approaches, (\ref{eqT1}) to (\ref{eqA2}) are substituted into one of the partial differential equations
\begin{align}
		\opcurl (\rho \opcurl \widetilde{\boldsymbol{T}}) + j \omega \mu \widetilde{\boldsymbol{T}} &= \boldsymbol{0} \text{ or} \label{eqA} \\
		\opcurl (\nu \opcurl \widetilde{\boldsymbol{A}}) + j \omega \sigma \widetilde{\boldsymbol{A}} &= \boldsymbol{J}_0 \label{eqT} 
\end{align}
and known steps considering the BCs (\ref{TMS1_BC_1}) to (\ref{AMS2_BC_3}) are carried out for the FEM (\hspace{-0.1cm} \cite{SchoeSchoeHoll_2D1D:19, HollSchoebi2D1DEE:20, HollSchoe:18}). \\ 
Finite element subspaces of the potentials have been selected as follows: $\boldsymbol{T}_0 \in H(\opcurl,\Omega)$, $\boldsymbol{T}_2$ and $\boldsymbol{A}_1 \in H(\opcurl,\Omega_m)$, $\Phi_0$ and $u_{10} \in H^1(\Omega)$ and $u_1$ and $w_{1} \in H^1(\Omega_m)$, see \cite{SchoeZagl:05}. The MSFs are in the space of periodic and continuous functions $H_{per}([-p/2,p/2])$. \\
Following the usual designations, we call $H(\opcurl)$ conforming FEs edge elements and $H^1$ conforming FEs nodal elements. The FE order (FEO) refers to edge elements for $\boldsymbol{T}_{0}$, $\boldsymbol{T}_{2}$ and $\boldsymbol{A}_{1}$. The FEO of nodal elements for $\Phi_0$, $u_{10}$, $u_1$ and $w_1$ in (\ref{eqA1}) is one higher than that of the edge elements to be consistent with the de-Rham-complex. An exception is $w_1$ in (\ref{eqA2}). \textcolor{blue}{MSFEM approach (\ref{eqT1}) uses only edge elements and (\ref{eqA2}) only nodal elements.} \\
The weak forms have been derived as described for example in \cite{HollSchoe:18} and \cite{Holl:19} for a MSFEM based on a MVP and in \cite{HansSchoeHoll:22} for a MSFEM based on a CVP. 
Averaging of the coefficients in the bilinear forms has been carried out to exploit the advantage of the MSFEMs, see \cite{SchoeSchoeHoll_2D1D:19}.
The BSF $\boldsymbol{H}_{BS}$ is included into the approaches (\ref{eqT1}) and (\ref{eqT2}) based on a CVP. In case of approaches based on a MVP the known current density $\boldsymbol{J}_0$ is directly considered on the right hand side in the corresponding linear form of the MSFEM and integration by parts.
\section{Reference Solutions} \label{RefSols}
Reference solutions have been computed using the mixed formulations $\boldsymbol{A},V$-$\boldsymbol{A}$ and $\boldsymbol{T},\Phi$-$\Phi$ (\hspace{-0.1cm} \cite{Biro:99}), where $V$ is the electric scalar potential and $\Phi$ the magnetic scalar potential, and with 3D FE models with second order FEs of the entire segment in Fig.~\ref{MotorSegmentGeo}. \\
Only mixed formulations, e.g. $\boldsymbol{T},\Phi$-$\Phi$ and $\boldsymbol{A},V$-$\boldsymbol{A}$, allow the modeling of all BCs of the considered specific problem. For example, the representation of BC on $\Gamma_o$ is not possible using only $\boldsymbol{A}$.
\subsubsection{Boundary Value Problem with $\bf{A}$,$V$-$\bf{A}$} \label{3D_AVA_BVPs}
\begin{eqnarray}
	\opcurl \nu \opcurl\boldsymbol{A}+j\omega\sigma (\boldsymbol{A}+\opgrad V) = \boldsymbol{0}&\text{and } \nonumber \\
	\opdiv (j\omega \sigma (\boldsymbol{A}+\opgrad V)) = 0&\text{on } \Omega_c \label{eq_AVA_1} \\
	\opcurl \nu \opcurl\boldsymbol{A} = \boldsymbol{J}_0&\text{on } \Omega_0 \label{eq_AVA_2}
\end{eqnarray}
\begin{eqnarray}
	\nu \opcurl \boldsymbol{A} \times \boldsymbol{n} = \boldsymbol{0}&\text{on } \Gamma_{H} \label{eq_AVA_3} \\
	-j\omega\sigma (\boldsymbol{A}+\opgrad V)\cdot \boldsymbol{n} = 0&\text{on } \Gamma_{J} \label{eq_AVA_4} \\
\boldsymbol{A} \times \boldsymbol{n} = \boldsymbol{0}&\text{on } \Gamma_{B} \label{eq_AVA_5} \\
	V = 0&\text{and } \nonumber \\
\boldsymbol{A} \times \boldsymbol{n} = \boldsymbol{0}&\text{on } \Gamma_{E} \label{eq_AVA_6}
	\end{eqnarray}
\subsubsection{Boundary Value Problem with $\bf{T}$,$\Phi$-$\Phi$} \label{3D_TPhiPhi_BVPs}
\begin{eqnarray}
	\opcurl \rho \opcurl\boldsymbol{T}+j\omega\mu (\boldsymbol{T}-\opgrad \Phi) &  \nonumber \\
	=-\opcurl \rho \opcurl\boldsymbol{H}_{BS}-j\omega\mu \boldsymbol{H}_{BS}&\!\!\!\text{and }  \nonumber \\
	j\omega \opdiv (\mu (\boldsymbol{T}-\opgrad \Phi)) = -j\omega \opdiv (\mu \boldsymbol{H}_{BS})&\!\!\!\text{on } \!\Omega_c \label{eq_TPhi_1} \\
	-j\omega \opdiv (\mu \opgrad \Phi) = -j\omega \opdiv (\mu \boldsymbol{H}_{BS})&\!\!\!\text{on } \!\Omega_0 \label{eq_TPhi_2}
\end{eqnarray}
\begin{eqnarray}
\boldsymbol{T} \times \boldsymbol{n} = \boldsymbol{0}&  \nonumber \\
	\Phi = 0&\text{on } \Gamma_{H} \label{eq_TPhi_3} \\
	\boldsymbol{T} \times \boldsymbol{n} = 0&\text{on } \Gamma_{J} \label{eq_TPhi_4}  \\
	-j\omega\mu (\boldsymbol{T}+\boldsymbol{H}_{BS}-\opgrad \Phi) \cdot \boldsymbol{n} = 0&\text{on } \Gamma_{B} \label{eq_TPhi_5} \\
	\rho \opcurl\boldsymbol{T} \times \boldsymbol{n} =  \boldsymbol{0}&\text{on } \Gamma_{E} \label{eq_TPhi_6}
		\end{eqnarray}
Details of the associated weak forms can be found, for example, in \cite{Biro:99}.
\section{Numerical Simulations} \label{NumSimulations}
\subsection{Problem} \label{NumProblem}
The problem is a segment of a fictitious electrical machine shown with details in Fig.~\ref{MotorSegmentGeo}. A conductivity of $\sigma=2.08\cdot10^6\mathrm{S/m}$, a relative permeability of $\mu_r=1,000$, a frequency of $f=50\mathrm{Hz}$, a thickness of the sheet of $d=0.5\mathrm{mm}$, a fill factor of $k_f=0.95$ and a peak value of the current $I=100\mathrm{A}$ have been selected. The BSF fulfills the rotational symmetry due to the prescribed currents (\ref{EQ_5}), half of which point into the opposite direction. \\
\subsection{Results} \label{NumResults}
\begin{figure}
		\begin{minipage}{0.99\columnwidth}
			\begin{minipage}{0.32\columnwidth}
					\includegraphics[scale=0.301]{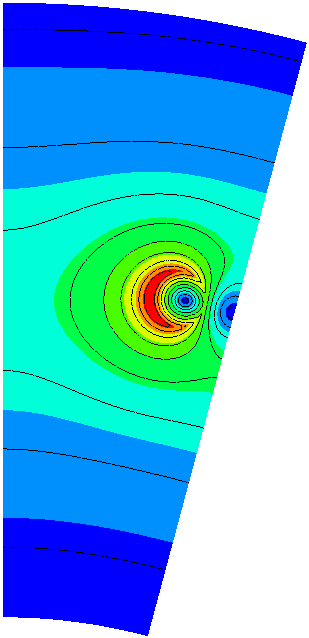} 
			\end{minipage}	
			\begin{minipage}{0.32\columnwidth}
					\includegraphics[scale=0.301]{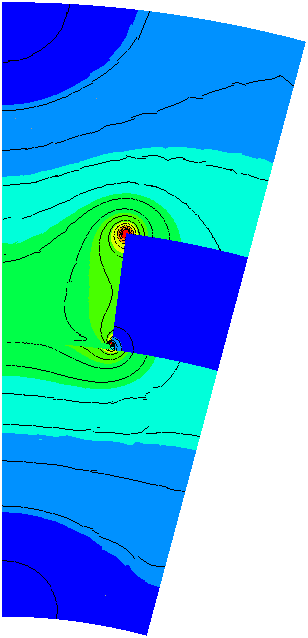} 
			\end{minipage}	
			\begin{minipage}{0.32\columnwidth}
					\includegraphics[scale=0.301]{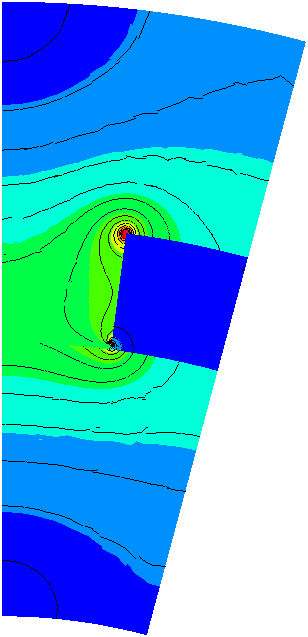} 
			\end{minipage}	
		\end{minipage}
	\vspace{0.3cm} \\
		\begin{minipage}{0.99\columnwidth}
			\begin{minipage}{0.32\columnwidth}
					\includegraphics[scale=0.301]{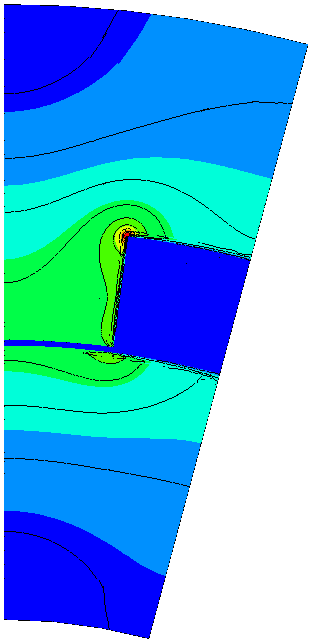} 
			\end{minipage}	
			\begin{minipage}{0.32\columnwidth}
					\includegraphics[scale=0.301]{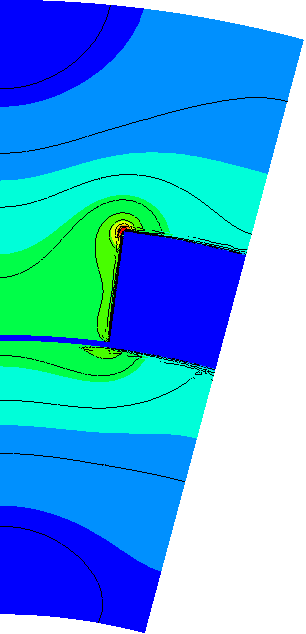} 
			\end{minipage}	
			\begin{minipage}{0.32\columnwidth}
					\includegraphics[scale=0.301]{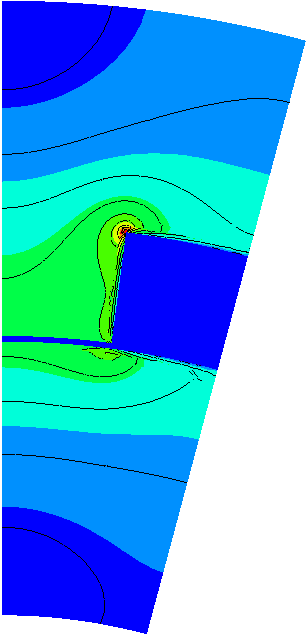} 
			\end{minipage}	
		\end{minipage}
\caption{Top: Biot-Savart-field $|\boldsymbol{H}_{BS}|$ to the left ($|\boldsymbol{H}_{BS}|_{max}=\mathrm{3,779kA/m}$), magnetic flux density $|\boldsymbol{B}|$ of T3D in the middle and of TMS1 to the right ($|\boldsymbol{B}|_{max}=\mathrm{0.273T}$), all at z=0. Bottom: Current density $|\boldsymbol{J}|$ ($|\boldsymbol{J}|_{max}=\mathrm{35kA/m^2}$) of A3D to the left, of AMS1 in the middle and of AMS2 to the right, all at z=0.188mm.}
\label{FieldDistribution}
\end{figure}

\begin{figure}
\begin{center} \includegraphics[scale=0.1825]{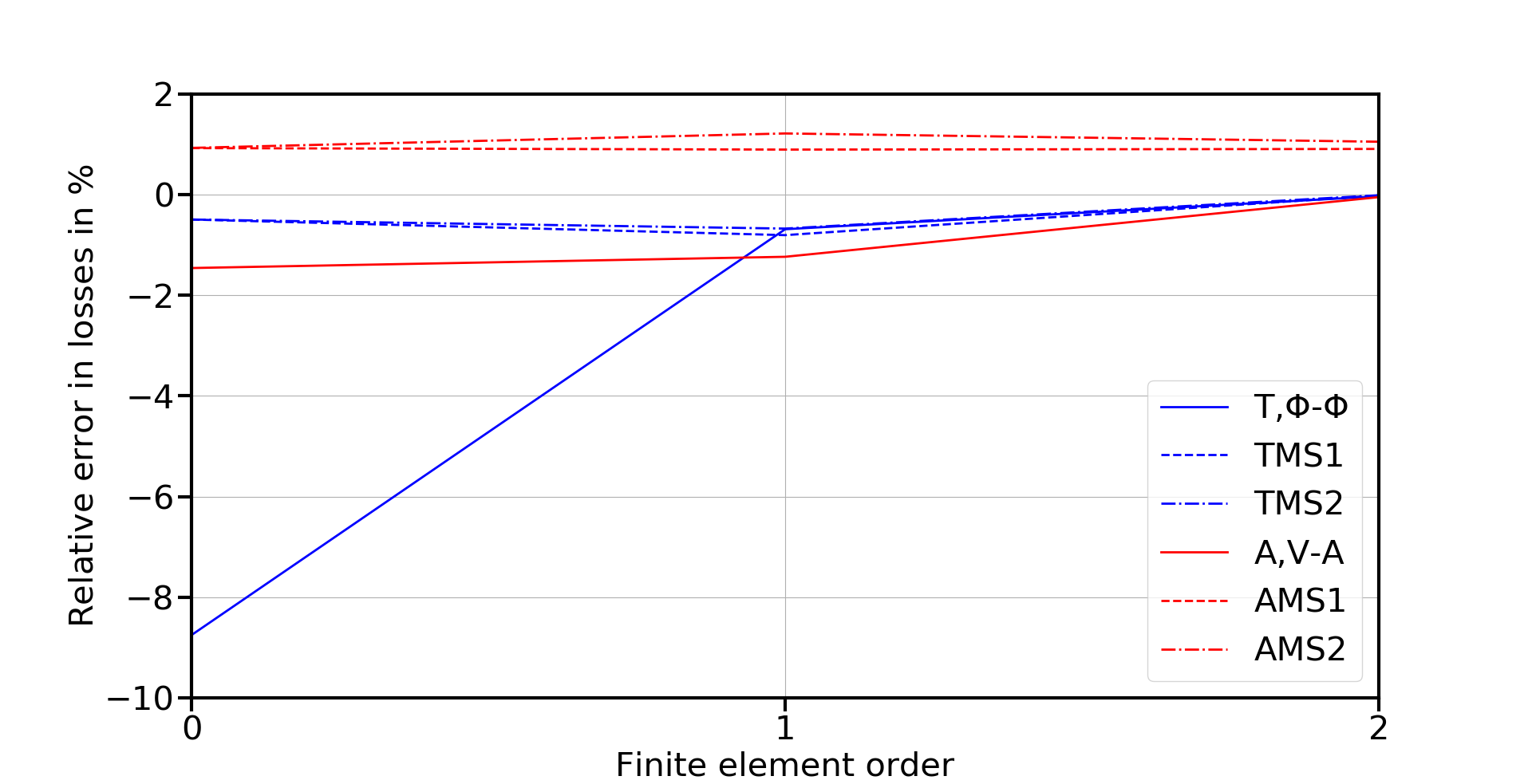} \end{center}
\caption{Relative error of eddy current losses.}
\label{RelErrorEddyLossestothalf}
\vspace{-0.105cm}
\end{figure}

Abbreviations TMS1, TMS2, AMS1 and AMS2 are used for solutions of MSFEMs based on the approaches (\ref{eqT1}), (\ref{eqT2}), (\ref{eqA1}) and (\ref{eqA2}), respectively. Reference solutions have been computed using the mixed formulations $\boldsymbol{A},V$-$\boldsymbol{A}$ and $\boldsymbol{T},\Phi$-$\Phi$ (\hspace{-0.1cm} \cite{Biro:99}), where $V$ is the electric scalar potential and $\Phi$ the magnetic scalar potential, and with 3D FE models with second order FEs of the entire segment in Fig.~\ref{MotorSegmentGeo}. Results obtained by TMS1 and TMS2 are evaluated by $\boldsymbol{T},\Phi$-$\Phi$ and those by AMS1 and AMS2 by $\boldsymbol{A},V$-$\boldsymbol{A}$. \\ 
For the sake of a fair comparison the  FE mesh for the 3D FEM has been generated by extrusion of the 2D mesh for the 2D/1D MSFEMs with prism-shaped elements. The 3D FEM models consist of six layers of elements with proper thicknesses to account for the penetration depth. \\ 
For comparison, some solutions are shown in Fig.~\ref{FieldDistribution}. There is a very satisfactory agreement. All methods have been implemented in Netgen/NGSolve \cite{netgenngsolve} and all problems have been solved by the direct solver PARDISO \cite{Schenk2011}. \\
The relative errors (REs) of the eddy current losses presented in Fig.~\ref{RelErrorEddyLossestothalf} show a very satisfactory agreement. Losses obtained by 2D/1D MSFEMs differ from the reference losses by less than one percent. The REs have been obtained for the half segment shown in Fig.~\ref{MotorSegmentGeo}. Errors of T3D and A3D are due to the exploitation of symmetry. The absolute value of the RE is well below 2$\%$, except for zero order $\boldsymbol{T},\Phi$-$\Phi$. 
The relatively large error for lowest order $\boldsymbol{T},\Phi$-$\Phi$ can be explained by the fact that the current density $\boldsymbol{J}$ is the circulation of the CVP $\opcurl \boldsymbol{T}$ being just piecewise constant with lowest order edge elements. The high accuracy of 2D/1D MSFEMs with zero order FEs is due to the local description of the solution using MSFs. The RE is negligible small for TMS1 and TMS2 and second order FEs. However, the RE is about 1$\%$ for AMS1 and AMS2 independent of the FEO. Simply speaking, this indicates that the formulations of AMS1 and AMS2 are less suitable. \\
For reference, the eddy current losses are summarized in Tab.~\ref{tab_losses}.
\begin{table}[b]
	\begin{center}
		\caption{Eddy Current Losses}
		\begin{tabular}{c|c|c}
			Method, Approach & $\boldsymbol{T},\Phi$-$\Phi$ & $\boldsymbol{A},V$-$\boldsymbol{A}$ \\ \hline
			$P_{}$ in $\mu$W & 47.65 & 47.40 \\ \hline
		\end{tabular} 
		\label{tab_losses} \\ 
	\end{center}
	Second order FEs \\
	3D FE models \\ 
\end{table}
Modeling of the EE requires just a homogenous BC for $\boldsymbol{T}_2$ for $\boldsymbol{T}$ formulations, which is exact for the MSFEM approaches (\ref{eqT1}) and (\ref{eqT2}), respectively. However, an additional term for $\boldsymbol{A}$ formulations in the MSFEM approaches is required, the third term in (\ref{eqA1}) and (\ref{eqA2}), respectively, which is an approximation only. \\ 
The 2D/1D approaches ignore any magnetic stray fields, except (\ref{eqT1}) could consider an axial magnetic flux if useful, compare with (\ref{eq2D1Daux_1}). The others miss this capability. \\
 
\begin{figure}
\begin{center} \includegraphics[scale=0.1825]{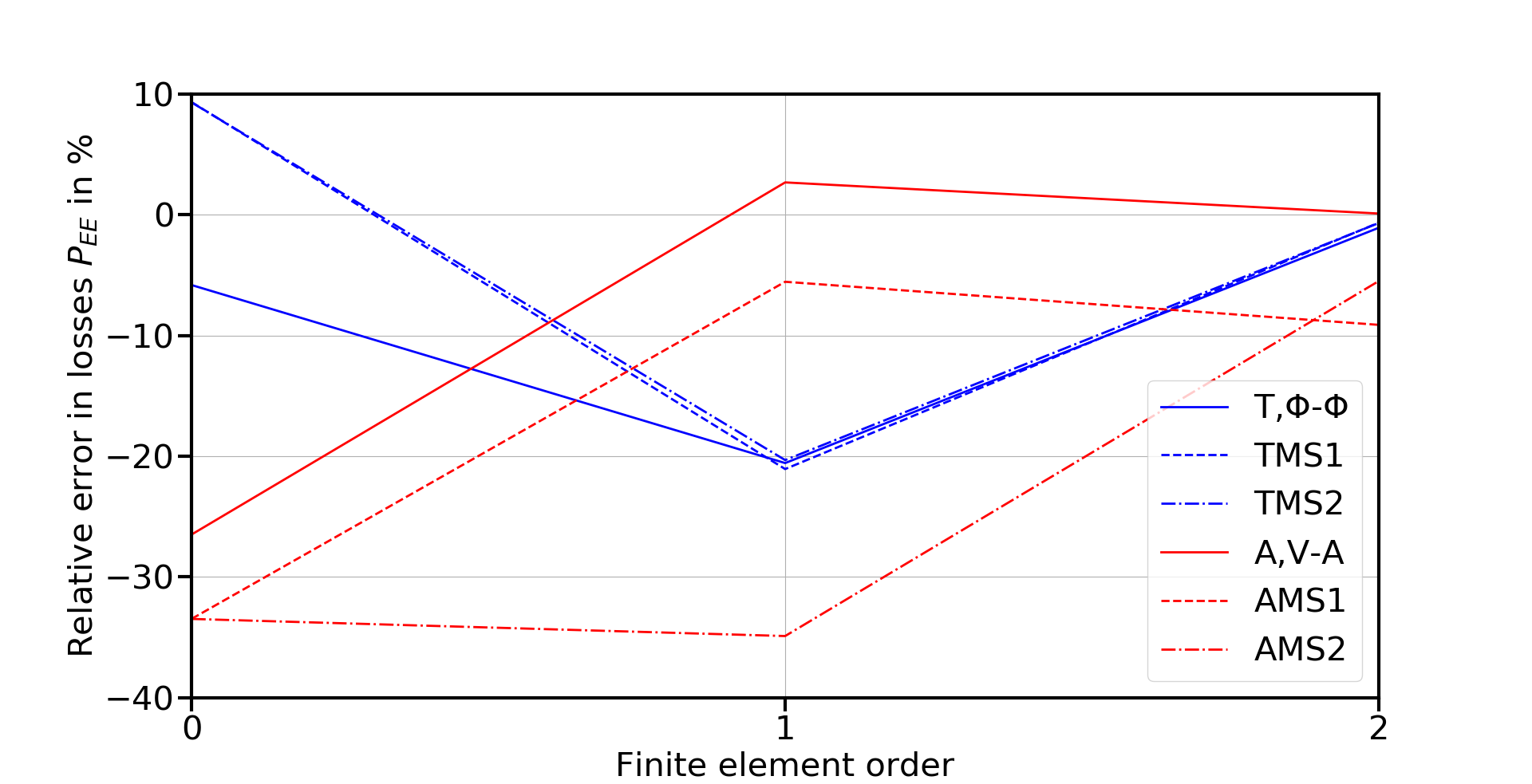} \end{center}
\caption{Relative error of losses $P_{EE}$ according to \ref{eqPEE}.}
\label{RelErrorEddyLossesEEhalf}
\vspace{-0.105cm}
\end{figure}

The main magnetic flux, which is parallel to the iron sheet, is considered by all approaches and causes eddy currents confined to flow in very narrow loops. The normal component, i.e. the $z$-component, of the current density represents the so-called EE. The EE is almost not present in the plane of symmetry. 
The losses
\begin{equation}
P_{EE} = 0.5 \int_{\Omega_c}{\sigma^{-1} J_z J_z^*~d\Omega},
\label{eqPEE}
\end{equation}
where $J_z$ is the $z$-component of the current density and * means conjugate complex, have been computed as a measure to study the capability to consider the EE by the 2D/1D MSFEMs. 
The RE of the EE practically vanishes for TMS1 and TMS2 with second order FEs, whereas it stays relatively large for AMS1 and AMS2, as can be seen in Fig.~\ref{RelErrorEddyLossesEEhalf}. \\
\begin{figure}
\begin{center} \includegraphics[scale=0.1825]{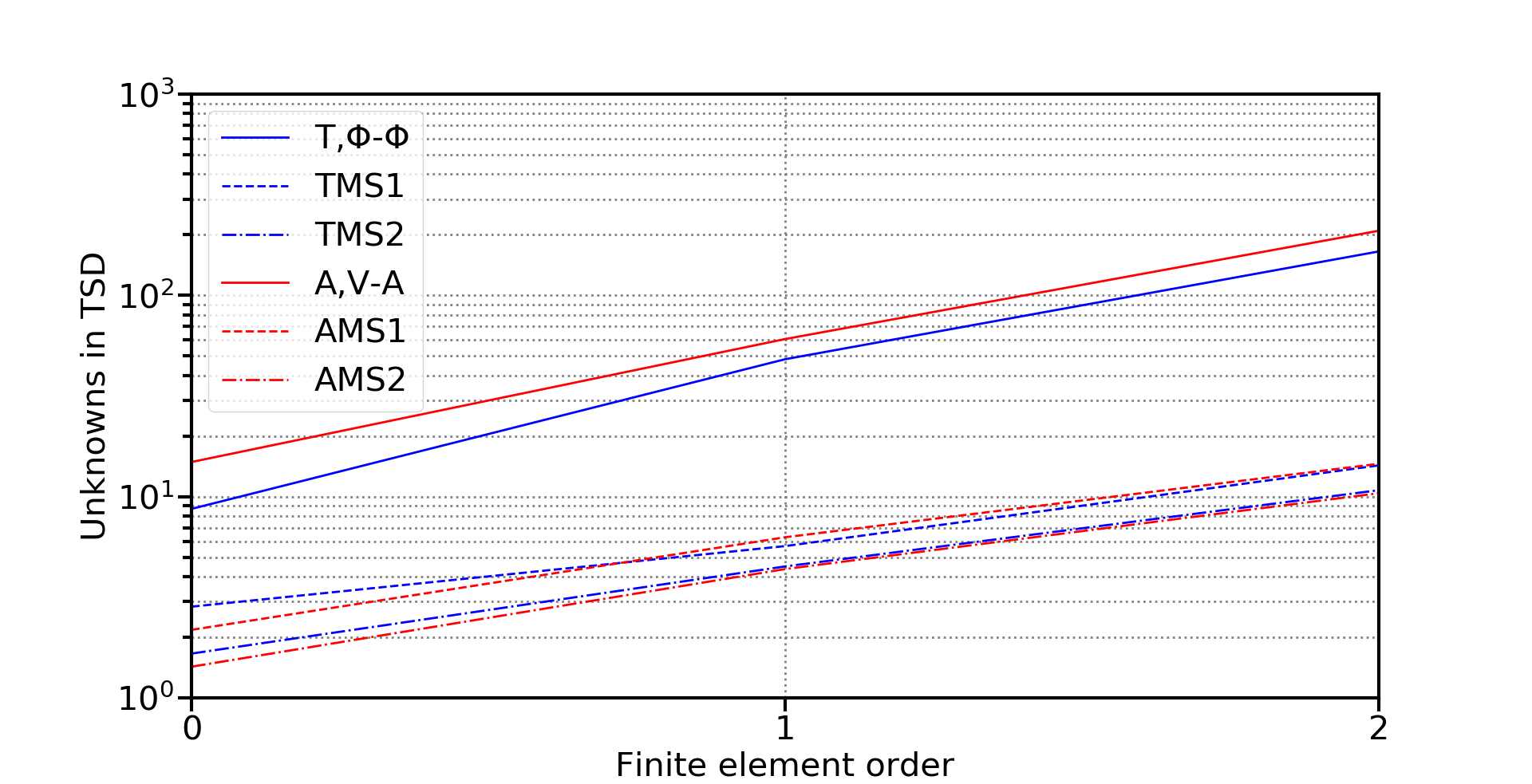} \end{center}
\caption{Number of unknows in thousands (TSD).}
\label{NoDOFshalf}
\vspace{-0.105cm}
\end{figure}
\begin{figure}
\begin{center} \includegraphics[scale=0.1825]{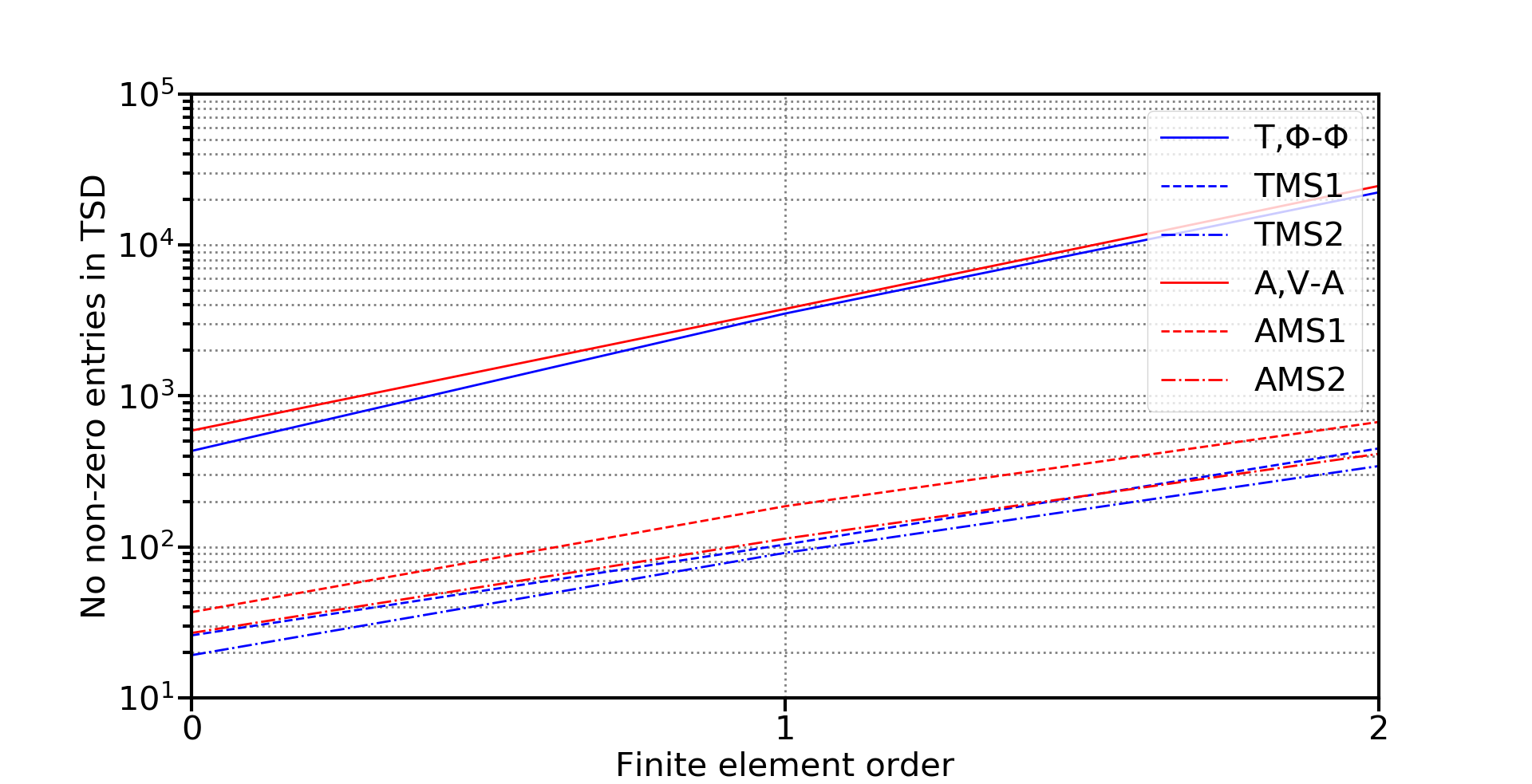} \end{center}
\caption{Non-zero entries in finite element matrix.}
\label{noZeroshalf}
\vspace{-0.105cm}
\end{figure}
\begin{figure}
\begin{center} \includegraphics[scale=0.1825]{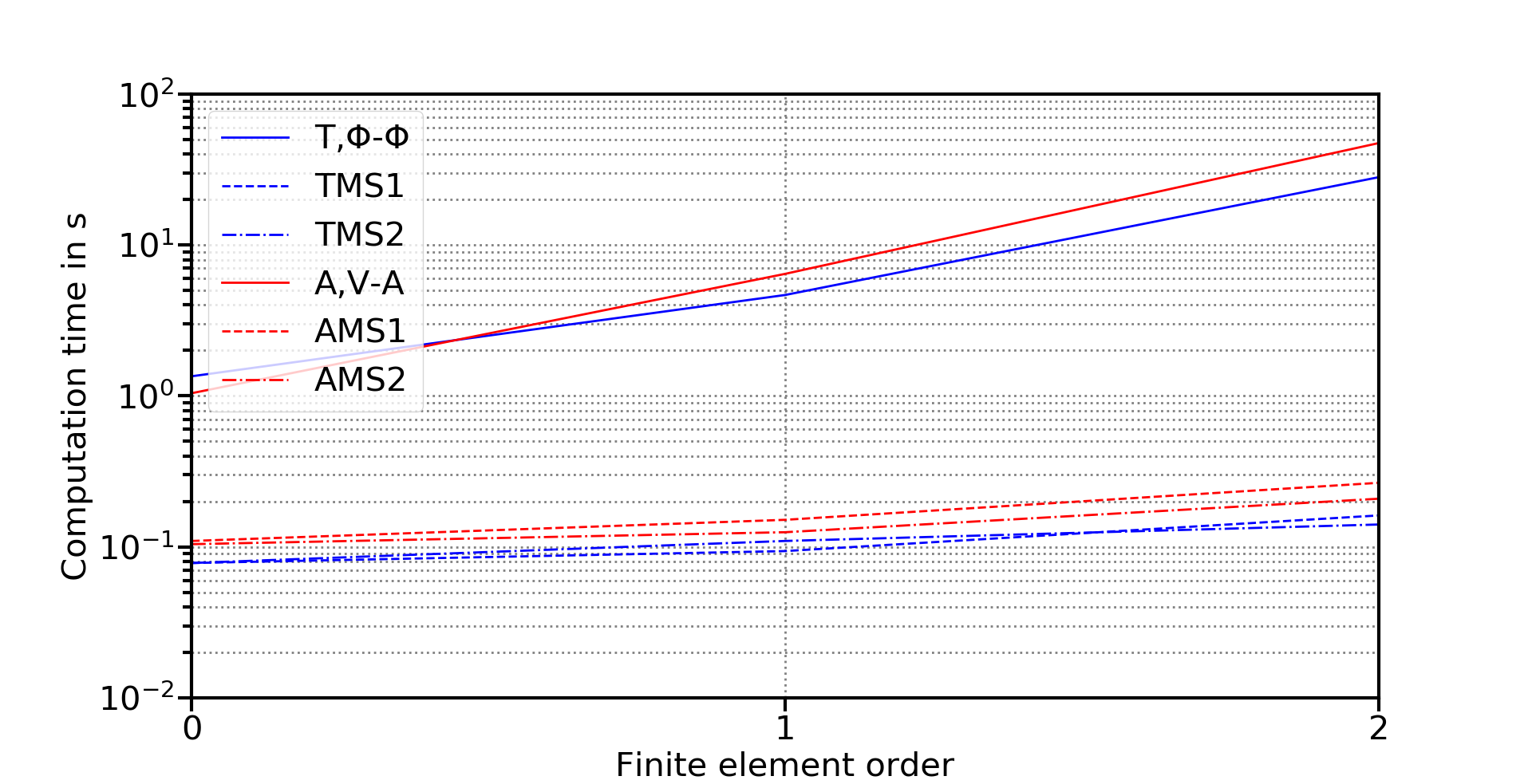} \end{center}
\caption{Computation time, half problem.}
\label{ComputTime_half}
\vspace{-0.105cm}
\end{figure}
\begin{figure}
\begin{center} \includegraphics[scale=0.1825]{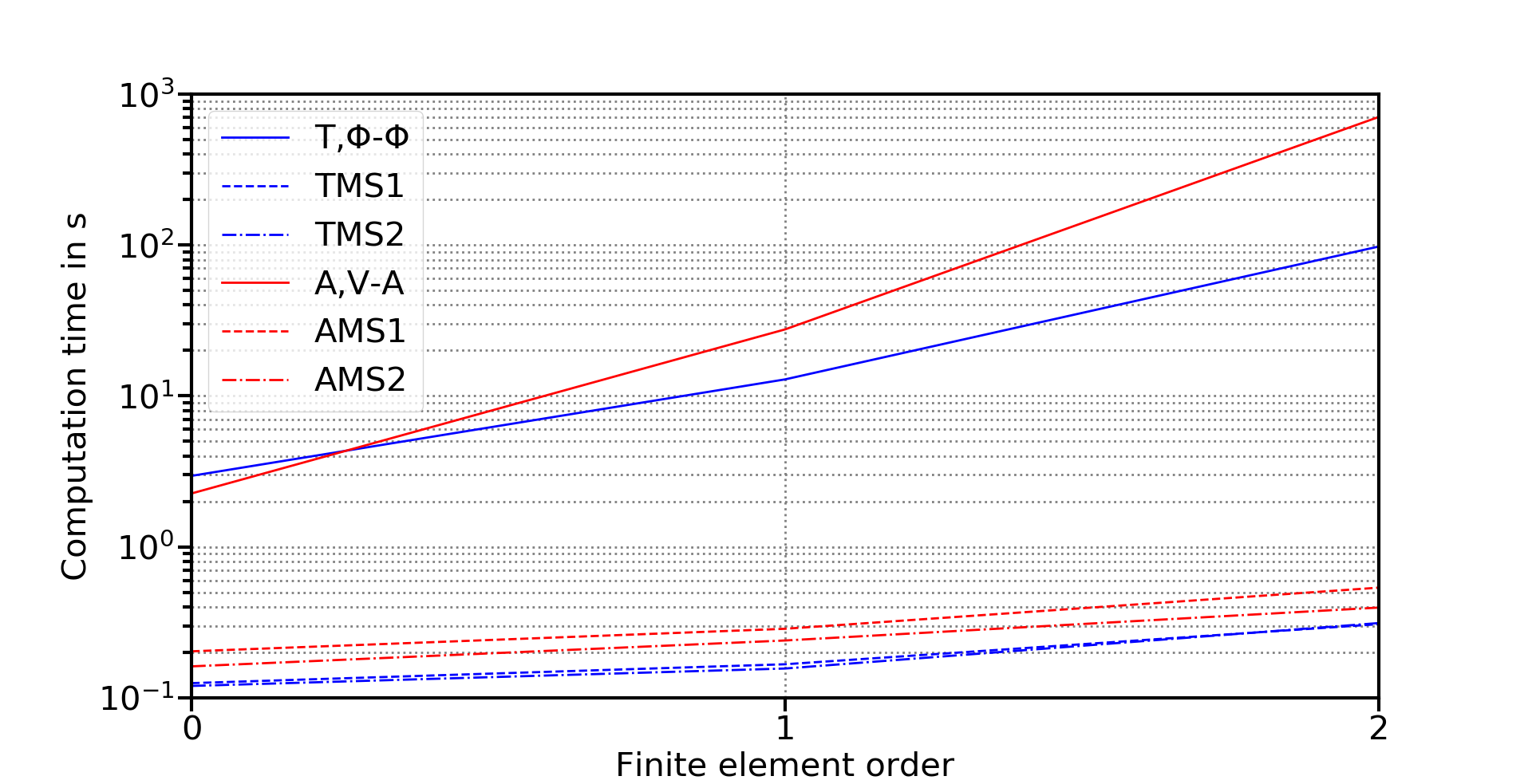} \end{center}
\caption{Computation time, entire problem.}
\label{ComputTime_entire}
\vspace{-0.105cm}
\end{figure}
The number of unknowns are presented in Fig.~\ref{NoDOFshalf}. In general, the 2D/1D MSFEMs require essentially less unknowns than the 3D FEMs, about a factor of 10. The $0^{th}$, the $1^{st}$ and the $2^{nd}$ order approximation belongs to the whole MSFEM approaches in (\ref{eqT1}) to (\ref{eqA2}). For example, $1^{st}$ order means for (\ref{eqA1}) that for $u_{10}$ $2^{nd}$ order, for $\boldsymbol{A}_1$ $1^{st}$ order and for $w_1$ $2^{nd}$ order standard finite element spaces are used. A significant additional reduction of the unknowns by the 2D/1D MSFEMs can be obtained by replacing $\boldsymbol{T}_{0}$ and $\boldsymbol{A}_{1}$ by $\opgrad \Phi_0$ and $\opgrad u_{1}$, respectively. To be fair the 3D FEMs could exploit the symmetry with respect to the plane z=0.

The memory requirement is reflected by means of the non-zero entries in Fig.~\ref{noZeroshalf}. Memory requirements increase visibly less rapidly for 2D/1D MSFEMs and are at least 10 times smaller than for 3D FEMs.

The computation times (CTs) presented in Figs.~\ref{ComputTime_half} and \ref{ComputTime_entire} consist of the solution of the problem including the evaluation of eddy current losses $P$ and $P_{EE}$. The CTs and their increase from the half to the entire problem with respect to the unknowns of the 3D reference solutions are essentially higher than those of the 2D/1D-MSFEM problems. Overall, the 2D/1D MSFEMs are more than 100 times faster than the 3D FEMs for second order FEO.
\section{Conclusions} \label{Conclusions}
The approaches for 2D/1D MSFEMs are quite similar to those for 3D MSFEMs \cite{Holl:19, HollSchoeb_TPhi:20}. While the coefficient functions for 2D/1D MSFEMs depend on two variables, those for 3D depend on three variables. Note, the term $\phi_1^0 \opgrad u_{10}$ in  (\ref{eqA1}) and (\ref{eqA2}) is not suitable in 3D. \\
Losses obtained by 2D/1D MSFEMs differ from the reference losses by less than about one percent. The MSFEMs with $\boldsymbol{T}$ formulations are noticeably more accurate than MSFEMs with $\boldsymbol{A}$ formulations. Overall, the number of unknowns for the 2D/1D MSFEMs is much smaller than that for the 3D FEMs. Simulations with 2D/1D MSFEMs are much faster than those with 3D FEMs. \\
The 2D/1D MSFEMs are able to handle problems with complicated geometries, Biot-Savart fields, symmetries and the EE. \\
Application of the proposed methods to the time domain or to nonlinear problems is obviously possible. \\
Therefore, using 2D/1D MSFEMs is a very attractive alternative to brute force 3D FEMs. 

\section*{Acknowledgment}
This work was supported by the Austrian Science Fund (FWF) under projects P 31926.
\bibliographystyle{IEEEtran}
\bibliography{hollaus}

\end{document}